# ON CHEBYSHEV SYSTEMS AND NON-UNIFORM SAMPLING RELATED TO CONTROLLABILITY AND OBSERVABILITY OF CAPUTO FRACTIONAL DIFFERENTIAL SYSTEMS

by M. DE LA SEN

Institute for Research and Development of Processes. Faculty of Science and Technology

University of Basque Country. Campus of Leioa. Aptdo. 544- Bilbao – SPAIN

**Abstract**. This paper is concerned with the investigation of the controllability and observability of Caputo fractional differential linear systems of any real order $\alpha$. Expressions for the expansions of the evolution operators in powers of the matrix of dynamics are first obtained. Sets of linearly independent continuous functions or matrix functions ,which are also Chebyshev systems, appear in such expansions in a natural way. Based on the properties of such functions, the controllability and observability of the Caputo fractional differential system of real order $\alpha$ are discussed as related to their counterpart properties in the corresponding standard system defined for $\alpha = 1$. Extensions are given to the fulfilment of those properties under non- uniform sampling. It is proved that the choice of the appropriate sampling instants ion not restrictive as a result of the properties of the associate Chebyshev systems.

## 1. Introduction

Caputo fractional calculus is a very useful tool to calculate alternative solutions to the classical ones in many applications as, for instance, in dynamic systems (see, for instance, [1-3]). Since the fractional order $\alpha$ can be a non positive integer, even real or complex, the technique can be used to better fix the trajectory solution of mathematical models to obtained experimental data due to modelling or measuring errors. On the other hand, the so- called Chebyshev system of n linearly independent functions in the Banach space of continuous functions $C([a,b])$, endowed with the supremum norm, has the important property, due to Haar, that each non- trivial polynomial polynomial of this system has at most (n-1) distinct zeros in $[a,b]$, [4]. There are many sets of linearly independent functions which are Chebyshev systems as , for instance, (a) The sets $\{t^k : k = 0,...,n-1\}$, respectively, $\{cos\, kt : k = 0,...,n-1\}$ on any real interval $[a,b]$ of nonzero finite measure, respectively, on $[0, 2\pi]$ and (b) the set $\{f(t)t^k : k = 0,...,n-1\}$ with f(t) being continuous and with no zeros on $[a,b]$ is a Chebyshev system on $[a,b]$. The linearly independent functions defining the expansions of the $C_0$- semigroup $\{e^{At} : t \in \mathbf{R}_{0+}\}$, of differential generator A, in powers of A are a Chebyshev system on any interval $[\gamma, \gamma + \pi/\omega)$ where $\omega$ is the maximum eigenfrequency ; i.e. the maximum absolute imaginary part of any complex eigenvalues, if any, and otherwise $\omega = +\infty$ (see [5-7], firstly discussed in [6]). In [5], it is proven that the property holds even if A is a complex matrix which is not related by a similarity transformation to some real one, i.e. some potential complex eigenvalues may not have their complex conjugate counterparts as eigenvalues. In [6-7], the Haar property of Chebyshev systems is used to formulate the following properties:

(a) controllability to the origin ( roughly speaking, the ability to steer any initial condition to the equilibrium in any finite time by injecting some admissible control) under, in general, non-uniform (also referred to as aperiodic or non-periodic) sampling

(b) observability (the ability to calculate initial conditions from past values of a measured output trajectory) under, in general, non-uniform sampling

(c) identifiability (the ability to compute the values of the parameters from output measurements and eventually time- derivatives up till some order of the output trajectory) and model matching of linear and time- invariant systems also under, in general, non-uniform sampling.

It has been proven in the above papers that the Haar property ensures that the corresponding properties are transferred from the continuous- time dynamic system to its discrete-time counterpart for infinitely many choices of sets of sufficient large cardinal (exceeding a lower-bound being dependent on the degree of the minimal polynomial of A and the dimensions of the output and input spaces) of distinct sampling instants. The particular choices of the sampling instants may be done by some practical considerations as, for instance, to achieve a well- conditioned coefficient matrix of the resulting algebraic problem  related to controllability, observability, local identifiability etc. The study of the properties of controllability, reachability (related to controllability to any final state), observability and constructability (closely related to observability) under  non-uniform sampling has been extended more recently to positive linear continuous time systems, namely, those having nonnegative state and output trajectory for any given non- negative initial conditions and controls,  [8].

The paper is organized as follows. Section 2 is devoted to some preliminary results concerning the expansions of the evolution operators of the solution in powers of the matrix of dynamics of  a linear time- invariant dynamic system.  Sets of Chebyshev systems appear in a natural way in such expansions. The results are extended in Section 3  for Caputo fractional dynamic systems on any real order $\alpha$ . Those results are used in Sections 4-5 to discuss the properties of observability and controllability and their counterparts under non-uniform sampling. Many other  theoretical studies and applications have been performed. In [10], the reconstruction problem from non-uniform data is focused on. In [11], the general non-uniform sampling in a stochastic framework is investigated. The problem of sequential estimation under non-uniform sampling is studied in [13]. In [12, 14-15], some estimation properties and application under non-uniform sampling are described. The sampling efficiency in event-based sampling laws of signals is discussed in [16]. Finally, some filtering properties under non-uniformly sampled signals are investigated in [17].

**I.1**. Some  notation  used  is   $R_{0+} := R \cup \{0\}$;   $R_+ := \{z \in R_+ : z \geq 0\}$.  A  corresponding  notation $Z_{0+}$ , $Z_+$  are used for  the corresponding subsets of the integer set **Z**.

$\bar{n} := \{1, 2, \cdots, n\}$ ; $\forall n \in Z_+$ .

## 2. Preliminaries on the expansions of exp(At) and associate Chebyshev systems

This section is devoted to preliminary basic results about useful sets of linearly independent real functions and real matrix functions of time, which are also Chebyshev systems  [4-7]), and  used for the expansions of the evolution operators defining the solution of the differential system.

**Lemma 2.1.** Assume that $\mu$ is the degree of the minimal polynomial of $A \in C^{n \times n}$. Then $\{A^i \in C^{n \times n} : i \in \overline{\mu - 1} \cup \{0\}\}$ is a linearly independent set of matrices.

**Proof**: Proceed by contradiction by assuming that the set $\{A^i \in C^{n \times n} : i \in \overline{\mu - 1} \cup \{0\}\}$ is linearly dependent so that there is a non-identically zero set $\{\lambda_i \in C : i \in \overline{\mu - 1} \cup \{0\}\}$ such that $\sum_{i=1}^{\mu-1} \lambda_i A^i = 0$. Assume that $\lambda_j \neq 0$ for $j \in \overline{\mu - 1} \cup \{0\}$ so that

$$A^j = -\lambda_j^{-1}\left(\sum_{i(<j)=1}^{\mu-1} \lambda_i A^i + \sum_{i(>j)}^{\mu-1} \lambda_i A^i\right) \tag{2.1}$$

Choose $j = max\left(i \in \overline{\mu-1} \cup \{0\} : \lambda_k = 0; \forall k(>i) \in \overline{\mu-1} \cup \{0\}\right)$ which always exists since $\{\lambda_i \in C : i \in \overline{\mu-1} \cup \{0\}\}$ is a non identically zero set and $0 \leq j \leq \mu - 1$. Thus, $A^j = -\sum_{i(<j)=1}^{\mu-1} \lambda_j^{-1} \lambda_i A^i$ and some $0 \leq j \leq \mu - 1$ is then the degree of the minimal polynomial of A what contradicts that such a degree is $\mu$. □

**Lemma 2.2.** The following formula $A^p = \sum_{k=0}^{\mu-1} a_k(p) A^k$ holds with unique complex coefficients $a_k(p)$; $\forall k \in \overline{\mu-1} \cup \{0\}, \forall p (\in Z_+) \geq \mu$.

**Proof**: Note from Cayley-Hamilton theorem that $A^\mu = \sum_{k=0}^{\mu-1} a_k(\mu) A^k$, thus,

$$A^{\mu+1} = \sum_{k=0}^{\mu-1} a_k(\mu) A^{k+1} = \sum_{k=0}^{\mu-2} a_k(\mu) A^{k+1} + a_{\mu-1} A^\mu$$

$$= \sum_{k=1}^{\mu-1} a_{k-1}(\mu) A^k + a_{\mu-1}(\mu)\left(\sum_{k=0}^{\mu-1} a_k(\mu) A^k\right)$$

$$= \sum_{k=0}^{\mu-1} \left(a_{k-1}(\mu) + a_k(\mu) a_{\mu-1}(\mu)\right) A^k = \sum_{k=0}^{\mu-1} a_k(\mu+1) A^k \tag{2.2}$$

so that the above identity holds irrespective of A if and only if:

$$a_k(\mu+1) = a_{\mu-1}(\mu) a_k(\mu) + a_{k-1}(\mu) = \frac{a_0(\mu+1)}{a_0(\mu)} a_k(\mu) + a_{k-1}(\mu); \forall k \in \overline{\mu-1} \cup \{0\}, a_{-1}(\mu) = 0 \tag{2.3}$$

Proceeding recursively with the above formulas, one gets for all $Z_+ \ni p \geq 1$

$$A^{\mu+p} = \sum_{k=0}^{\mu-1} a_k(\mu+p-1) A^{k+1} = \sum_{k=0}^{\mu-1} \left(a_{k-1}(\mu+p-1) + a_k(\mu+p-1) a_{\mu-1}(\mu+p-i-1)\right) A^k$$

$$= \sum_{k=0}^{\mu-1} a_k(\mu+p) A^k \tag{2.4}$$

with

$$a_k(\mu+p) = a^p_{\mu-1}(\mu)a_k(\mu) + \sum_{i=0}^{p-1} a^i_{\mu-1}(\mu)a_{k-1}(\mu+p-i-1) \tag{2.5}$$

The uniqueness of the coefficients follows from Lemma 2.1 as follows. Assume that there are two sets of coefficients such that $A^p = \sum_{k=0}^{\mu-1} a_k(p)A^k = \sum_{k=0}^{\mu-1} a'_k(p)A^k$. Then:

$$\sum_{k=0}^{\mu-1} \left(a_k(p) - a'_k(p)\right)A^k = 0 \text{ which implies } a_k(p) = a'_k(p); \quad \forall k \in \overline{p-1} \cup \{0\}; \quad \forall p(\in Z_+) \geq \mu \text{ from}$$

Lemma 2.1 since $\{A^i \in C^{n \times n} : i \in \overline{\mu-1} \cup \{0\}\}$ is a linearly independent set.

**Alternative proof**: Equivalently, proceed by complete induction by assuming that

$$A^{\mu+j} = \sum_{k=0}^{\mu-1} a_k(\mu)A^{k+j} = \sum_{k=0}^{\mu-1} a_k(\mu+j-1)A^{k+1} \; ; \; \forall j \in \overline{p-1} \tag{2.6}$$

with

$$a_k(\mu+j) = a^j_{\mu-1}(\mu)a_k(\mu) + \sum_{i=0}^{j-1} a^i_{\mu-1}(\mu)a_{k-1}(\mu+j-i-1) \; ; \; \forall j \in \overline{p-1} \tag{2.7}$$

Then, it follows that (2.4) holds with

$$a_k(\mu+j) = a_{\mu-1}(\mu+j-1)a_k(\mu+j-1) + a_{k-1}(\mu+j-1)$$
$$= a^j_{\mu-1}(\mu)a_k(\mu) + \sum_{i=0}^{j-1} a^i_{\mu-1}(\mu)a_{k-1}(\mu+j-i-1); \; \forall k \in \overline{\mu-1} \cup \{0\}, \; a_{-1}(j) = 0, \; \forall j \in \overline{p-1}$$

$$(2.8) \quad \square$$

A particular case of interest of (2.8) is

$$a_k(n) = a_{\mu-1}(n-1)a_k(n-1) + a_{k-1}(n-1)$$
$$= a^{n-\mu}_{\mu-1}(\mu)a_k(\mu) + \sum_{i=0}^{n-\mu-1} a^i_{\mu-1}(\mu)a_{k-1}(n-i-1) ; \; \forall k \in \overline{\mu-1} \cup \{0\}, \; a_{-1}(n-1) = 0 \tag{2.9}$$

what allows to interpret the two following common formulas for the Cayley-Hamilton theorem:

$$A^\mu = \sum_{k=0}^{\mu-1} a_k(\mu)A^k \; ; \; A^n = \sum_{k=0}^{\mu-1} a_k(n)A^k \tag{2.10}$$

where $\mu \leq n$ and n are, respectively, the degrees of the minimal and characteristic polynomials of A. On the other hand, it is possible to extend the expansion of powers of A up till the degree of the characteristic polynomial of $A$ through the use of modified coefficients $\bar{a}_k(.)$ as follows:

$$A^n = \sum_{k=0}^{n-1} \bar{a}_k(n)A^k = \sum_{k=0}^{\mu-1} a_k(\mu)A^k$$
$$= \sum_{k=0}^{\mu-1} \bar{a}_k(n)A^k + \sum_{k=\mu}^{n-1} \bar{a}_k(n)A^k = \sum_{k=0}^{\mu-1} \bar{a}_k(n)A^k + \sum_{k=0}^{n-1-\mu} \bar{a}_{k+\mu}(n)A^{k+\mu}$$

$$= \sum_{k=0}^{\mu-1} \overline{a}_k(n) A^k + \sum_{k=0}^{n-1-\mu} \overline{a}_{k+\mu}(n) \left( \sum_{j=0}^{\mu-1} a_j(\mu+k) A^j \right)$$

$$= \sum_{j=0}^{\mu-1} \left[ \overline{a}_j(n) + \sum_{k=0}^{n-1-\mu} \overline{a}_{k+\mu}(n) a_j(\mu+k) \right] A^j \qquad (2.11)$$

so that $\overline{a}_j(n) = a_j(\mu) - \sum_{k=\mu}^{n-1} \overline{a}_k(n) a_j(k)$ ; $\forall k \in \overline{n-1} \cup \{0\}$ and it turns out that the extension also works for $A^p$; $\forall p (\in \mathbf{Z}_+) \geq \mu$ with $\overline{a}_k(\mu) = a_k(\mu)$; $\forall k \in \overline{p-1} \cup \{0\}$, and

$$\overline{a}_k(p) = a_k(\mu) - \sum_{j=\mu}^{p-1} \overline{a}_j(p) a_k(k) \; ; \; \forall k \in \overline{p-1} \cup \{0\} \qquad (2.12)$$

for $p \geq \mu + 1$. It has been proven the following:

**Lemma 2.3**. There exist complex coefficients $a_k(p)$ and $\overline{a}_k(p); \forall k \in \overline{p-1} \cup \{0\}$ ; $\forall p (\in \mathbf{Z}_+) \geq \mu$, such that the formulas $A^p = \sum_{k=0}^{\mu-1} a_k(p) A^k = \sum_{k=0}^{p-1} \overline{a}_k(p) A^k$ are true. Those sets of coefficients coincide if $p = \mu$. □

It turns out that he sets of expanding coefficients are real if the matrix A is real or if it is complex being similar to a real one, i.e. all complex values, if any, appear by complex conjugate pairs or identical multiplicities. Now, the above results are related to the fundamental matrices of time- invariant differential systems as follows. Consider $C_0$- semigroup $\{\Phi(t) := e^{At} : t \in \mathbf{R}_{0+}\}$ of infinitesimal generator $A \in \mathbf{C}^{n \times n}$ is an evolution operator in $L(\mathbf{C}^n)$ point-wise defined by $\Phi(t) = e^{At}$ the defining the trajectory solutions $x(t) = \Phi(t) x_0$ of the differential system of order n $\dot{x}(t) = A x(t); x(0) = x_0$. An equivalent description is that $\Phi: \mathbf{R}_{0+} \to \mathbf{C}^{n \times n}$ is the fundamental matrix function (or the state-transition matrix function) of the above differential system. The expansions of $\Phi(t)$ in finite powers of A, being non less than its minimal polynomial, where studied in detail for $A \in \mathbf{C}^{n \times n}$ in [4-7] by using Chebyshev sets of complex functions.

**Lemma 2.4.** The formulas $\Phi(t) := e^{At} = \sum_{k=0}^{p-1} \beta_k(t,p) A^k$; $\forall p \geq \mu$ for linearly independent sets of functions $B_p := \{\beta_k : \mathbf{R}_{0+} \times \overline{p} \to \mathbf{R} : k \in \overline{p-1} \cup \{0\}\}$ which are unique and analytic in $\mathbf{R}_+$ for each $p \geq \mu$ and which satisfy a p-order linear time-invariant differential system on $\mathbf{R}_+$.

**Proof**: It follows after using Lemma 2.3 with $A^p = \sum_{k=0}^{\mu-1} a_k(p) A^k = \sum_{k=0}^{p-1} \overline{a}_k(p) A^k$ that

$$\Phi(t) := e^{At} = \sum_{k=0}^{\infty} \frac{A^k t^k}{k!} = \sum_{k=0}^{p-1} \frac{A^k t^k}{k!} + \sum_{k=p}^{\infty} \frac{A^k t^k}{k!} = \sum_{k=0}^{\mu-1} \frac{A^k t^k}{k!} + \sum_{k=\mu}^{\infty} \frac{A^k t^k}{k!}$$

$$= \sum_{k=0}^{\mu-1} \left( \frac{t^k}{k!} + \sum_{i=\mu}^{\infty} \frac{a_k(i)t^i}{i!} \right) A^k = \sum_{k=0}^{\mu-1} \beta_k(t,\mu) A^k$$

$$= \sum_{k=0}^{p-1} \left( \frac{t^k}{k!} + \sum_{i=p}^{\infty} \frac{\overline{a}_k(i)t^i}{i!} \right) A^k = \sum_{k=0}^{p-1} \beta_k(t,p) A^k \ , \forall p \geq \mu \tag{2.13}$$

with $\beta_0(0,p)=1$, $\beta_i(0,p)=0$ ; $\forall i \in \overline{p-1} \cup \{0\}$, where

$$\beta_k(t,\mu) = \frac{t^k}{k!} + \sum_{i=\mu}^{\infty} \frac{a_k(i)t^i}{i!} \ ; \ k \in \overline{\mu-1} \cup \{0\} \tag{2.14}$$

$$\beta_k(t,p) = \frac{t^k}{k!} + \sum_{i=p}^{\infty} \frac{\overline{a}_k(i)t^i}{i!} \ ; \ k \in \overline{p-1} \cup \{0\}, \ \forall p \geq \mu \tag{2.15}$$

Using the recursive expressions (2.9) for the $a_{(\cdot)}(\cdot)$ in the above formula and using also the recursions (2.12) for the $\overline{a}_{(\cdot)}(\cdot)$; $\forall k \in \overline{p-1} \cup \{0\}$, one gets:

$$\beta_k(t,\mu) = \frac{t^k}{k!} + \sum_{i=\mu}^{\infty} \frac{t^i}{i!} \left( a_{k-\mu-1}^{i-k+\mu}(k-\mu) a_k(k-\mu) + \sum_{\ell=0}^{i-k+\mu-1} a_{k-\mu-1}^{\ell}(k-\mu) a_{k-1}(i-\ell-1) \right)$$

$$; \ k \in \overline{\mu-1} \cup \{0\}, a_{-1}(i)=0 \ ; \ \forall i \geq \mu \tag{2.16}$$

$$\beta_k(t,p) = \frac{t^k}{k!} + \sum_{i=p}^{\infty} \left( a_k(\mu) - \sum_{j=\mu}^{i-1} \overline{a}_j(i) a_k(k) \right) \frac{t^i}{i!} \ ; \ k \in \overline{\mu-1} \cup \{0\} \ , \forall p \geq \mu \tag{2.17}$$

These functions satisfy sets of linear time – invariant ordinary differential equations as follows:

$$\dot{\Phi}(t) = A\Phi(t) = \sum_{k=0}^{p-1} \dot{\beta}_k(t,p) A^k = \sum_{k=0}^{p-1} \beta_k(t,p) A^{k+1}$$

$$= \sum_{k=1}^{p} \beta_{k-1}(t,p) A^k + \beta_{p-1}(t,p) A^p \ ; \ p \geq \mu \tag{2.18}$$

what implies

$$\dot{\beta}_0(t,p) I_n + \sum_{k=1}^{p-1} \left( \dot{\beta}_k(t,p) - \beta_{k-1}(t,p) \right) A^k - \beta_{p-1}(t,p) A^p = 0 \tag{2.19}$$

the above constraint holds for any matrix A so that functions in the sets $\boldsymbol{B}_p := \{ \beta_\ell : \boldsymbol{R}_{0+} \times \overline{p} \to \boldsymbol{R} : \ell \in \overline{p-1} \cup \{0\} \}$, $\forall p \geq \mu$ are linearly independent since they are the solutions of the linear time-invariant differential system:

$$\begin{bmatrix} \dot{\beta}_0(t,p) \\ \dot{\beta}_1(t,p) \\ \vdots \\ \dot{\beta}_{p-1}(t,p) \end{bmatrix} = \begin{bmatrix} 0 & 0 & \cdots & 0 & \overline{a}_0(p) \\ 1 & 0 & \cdots & 0 & \overline{a}_1(p) \\ & \ddots & \ddots & & \vdots \\ 0 & \cdots & 0 & 1 & \overline{a}_{p-1}(p) \end{bmatrix} \begin{bmatrix} \beta_0(t,p) \\ \beta_1(t,p) \\ \vdots \\ \beta_{p-1}(t,p) \end{bmatrix} \ ; \ \beta_0(0,p)=1, \ \beta_i(0,p)=0 \ ; \forall i \in \overline{p-1} \cup \{0\}$$

$$\tag{2.20}$$

; $\forall p \geq \mu, \forall k \in \overline{p-1} \cup \{0\}$ with $\overline{a}_k(\mu) = a_k(\mu); \forall k \in \overline{\mu-1} \cup \{0\}$, whose respective unique solutions are:

$$f_\alpha(t,p) = e^{\Omega_a(p)t} f_\alpha(0,p) \qquad (2.21)$$

where $f_\beta(t,p) = (\beta_0(t,p), \beta_1(t,p), \cdots, \beta_{p-1}(t,p))^T$ with $f_\beta(0,p) = e_1 = (1,0,\cdots,0)^T$, and

$$\Omega_a(p) = \begin{bmatrix} 0 & 0 & \cdots & 0 & \bar{a}_0(p) \\ 1 & 0 & \cdots & 0 & \bar{a}_1(p) \\ & \ddots & \ddots & & \vdots \\ 0 & \cdots & 0 & 1 & \bar{a}_{p-1}(p) \end{bmatrix} \qquad (2.22)$$

Then, there exist unique sets $B_p := \{\beta_k : R_{0+} \times \overline{p} \to R : k \in \overline{p-1} \cup \{0\}\}$ of linearly independent functions, which are unique and analytic in $R_+$ for each $p \geq \mu$ from (2.21)-(2.22), such that the formulas $\Phi(t) = e^{At} = \sum_{k=0}^{p-1} \beta_k(t,p) A^k$ hold; $\forall t \in R_{0+}$, $\forall p \geq \mu$. □

**Remark 2.5**: Note that it has been proven that there exist sets $A_p := \{a_k(p) : k \in \overline{\mu-1} \cup \{0\}\}$ and $\overline{A}_p := \{\bar{a}_k(p) : k \in \overline{p-1} \cup \{0\}\}$; $p \geq \mu$ of complex coefficients, with $\bar{a}_k(\mu) = a_k(\mu); \forall k \in \overline{\mu-1} \cup \{0\}$, such that the following polynomial expansions are true $A^\mu = \sum_{k=0}^{\mu-1} a_k(\mu) A^k$, $A^n = \sum_{k=0}^{n-1} \bar{a}_k(n) A^k = \sum_{k=0}^{\mu-1} a_k(n) A^k$ where $\mu$ and n are the degrees of the minimal and the characteristic polynomial of A, respectively. In general, $A^p = \sum_{k=0}^{\mu-1} a_k(\mu) A^k = \sum_{k=0}^{\mu-1} a_k(p) A^k = \sum_{k=0}^{p-1} \bar{a}_k(p) A^k$ for $p \geq \mu$. Those sets are unique for each $p \geq \mu$. The above formulas imply that there exist unique sets $B_p := \{\beta_k : R_{0+} \times \overline{p} \to R : k \in \overline{p-1} \cup \{0\}\}$ of linearly independent functions, which are unique and analytic in $R_+$ for each $p \geq \mu$ from (2.21)-(2.22), such that the formulas $\Phi(t) = e^{At} = \sum_{k=0}^{p-1} \beta_k(t,p) A^k$ hold; $\forall p \geq \mu$. □

The following result is useful to relate the coefficients of the minimal and characteristic polynomials of A to the two more relevant versions of the Cayley-Hamilton theorem:

**Lemma 2.6.** Assume that $p_n(A) = det(sI_n - A) = det(sI_n - J_A)$ is the characteristic polynomial of a nonzero matrix A of n order and Jordan form $J_A$. Then, $p_n(A) = s^n - \sum_{k=0}^{n-1} \bar{a}_k(n) s^k$.

Assume that $\mu$, with $1 \leq \mu \leq n$, is the degree of the minimal polynomial of A, $p_\mu(A)$. Then,

$p_\mu(A) = p_\mu(J_A) = p_\mu(\hat{J}_A) = s^\mu - \sum_{k=0}^{\mu-1} a_k(\mu) s^k$ where $\hat{J}_A = E^T J_A E = E^T T^{-1} J_A T E$ is a $\mu$-square matrix and $T$ and $E$ are n-square non-singular and $\mu \times n$ matrices, respectively, which are unique except multiplication by a nonzero scalar and $p_\mu(A)$ is the characteristic and minimal polynomial of $\hat{J}_A$.

**Proof**: The first part follows from the Cayley-Hamilton theorem in the form $A^n = \sum_{k=0}^{n-1} \bar{a}_k(n) A^k$ and the fact that the matrix satisfies its own characteristic polynomial. If $\mu$ is the degree of the minimal polynomial of A then there exist a non-singular real n-matrix T and a $\mu \times n$ real matrix E, formed by identity and zero block matrices of appropriate orders which pickup from $J_A$ the higher-order Jordan blocks of each eigenvalue of A, whose orders are given by the respective index of each eigenvalue, such that the identities $\hat{J}_A = E^T J_A E = E^T T^{-1} J_A T E$ hold. The last identity holds from the Cayley-Hamilton theorem since the characteristic polynomial of $\hat{J}_A$ and the minimal polynomial of A and $J_A$ coincide. □

### 3. The expansions of evolution operators of functional fractional Caputo differential systems and their associate Chebyshev systems

The above results are extended to the following fractional Caputo differential systems of order $\alpha$:

$$\left({}^C D_{0+}^\alpha x\right)(t) = Ax(t) + Bu(t) \tag{3.1}$$

with $k-1 < \alpha (\in \mathbf{R}_+) \leq k$; for some $k-1, k \in \mathbf{Z}_{0+}$ and $B \in \mathbf{R}^{n \times m}$ is the control matrix. If $\alpha = 1$ then (3.1) is referred to in the sequel as the *standard system*. The initial conditions are $\varphi_j(0) = x_j(0) = x^{(j)}(0) = x_{0j}$; $\forall j \in \overline{k-1} \cup \{0\}$. The admissible function vector $u: \mathbf{R}_{0+} \to \mathbf{R}^m$ is any given bounded piecewise continuous control function. The following result is concerned with the unique solution on $\mathbf{R}_{0+}$ of the above differential fractional system (3.1). The proof follows directly from a parallel existing result from the background literature on fractional differential systems by grouping all the additive forcing terms of (3.1) in a unique one (see, for instance [1], Eqs. (1.8.17), (3.1.34)-(3.1.49), with $f(t) \equiv Ax(t) + Bu(t)$). The linear and time-invariant differential functional fractional differential system (3.1) of any order $\alpha \in \mathbf{C}_{0+}$ has the following unique solution on $\mathbf{R}_{0+}$ for each given set of initial conditions and each given control $u: \mathbf{R}_{0+} \to \mathbf{R}^m$ being a bounded piecewise continuous control function:

$$x_\alpha(t) = \sum_{j=0}^{k-1} t^j \Phi_{\alpha j}(t) x_{0j} + \int_0^t (t-\tau)^{\alpha-1} \Phi_\alpha(t-\tau) Bu(\tau) d\tau; \ t \in \mathbf{R}_{0+} \tag{3.2}$$

with $k = [Re\, \alpha] + 1$ if $\alpha \notin \mathbf{Z}_+$ and $k = \alpha$ if $\alpha \in \mathbf{Z}_+$. The matrix functions $\Phi_{\alpha j}(t)$, $\forall j \in \overline{k-1} \cup \{0\}$ and $\Phi_\alpha(t)$ from $\mathbf{R}_{0+}$ to $\mathbf{R}^{n \times n}$ are calculated via the Mittag-Leffler functions which, after using the identities $A^\ell = \sum_{\ell=0}^{\mu-1} a_\ell(\ell) A^\ell = \sum_{\ell=0}^{p-1} \bar{a}_\ell(\ell) A^\ell$; $\forall p \geq \mu$ from Lemma 2.3 for the matrix $A$, lead to:

$$\Phi_{\alpha j}(t) = \sum_{\ell=0}^\infty \frac{A_0^\ell t^{\alpha \ell}}{\Gamma(\alpha \ell + j + 1)} = \sum_{\ell=0}^{\mu-1} \left( \frac{t^{\alpha \ell}}{\Gamma(\alpha \ell + j + 1)} + \sum_{i=\mu}^\infty \frac{\bar{a}_\ell(i) t^{\alpha i}}{\Gamma(\alpha i + j + 1)} \right) A^\ell$$

$$= \sum_{\ell=0}^{\mu-1} \beta_{\alpha j\ell}(t,\mu) A^{\ell} \ ; \ \forall j \in \overline{k-1} \cup \{0\}$$

$$= \sum_{\ell=0}^{p-1} \left( \frac{t^{\alpha\ell}}{\Gamma(\alpha\ell+j+1)} + \sum_{i=p}^{\infty} \frac{\overline{a}_{\ell}(i) t^{\alpha i}}{\Gamma(\alpha i+j+1)} \right) A^{\ell}$$

$$= \sum_{\ell=0}^{p-1} \beta_{\alpha j\ell}(t,p) A^{\ell} \ ; \ \forall p \geq \mu, \ \forall j \in \overline{k-1} \cup \{0\} \qquad (3.3)$$

$$\Phi_{\alpha}(t) = \sum_{\ell=0}^{\infty} \frac{A^{\ell} t^{\alpha\ell}}{\Gamma((\alpha+1)\ell)} = \sum_{\ell=0}^{p-1} \frac{A^{\ell} t^{\alpha\ell}}{\Gamma((\alpha+1)\ell)} + \sum_{\ell=p}^{\infty} \frac{A^{\ell} t^{\alpha\ell}}{\Gamma((\alpha+1)\ell)}$$

$$= \sum_{\ell=0}^{p-1} \left( \frac{t^{\alpha\ell}}{\Gamma((\alpha+1)\ell)} + \sum_{i=p}^{\infty} \frac{\overline{a}_{k}(\ell) t^{\alpha i}}{\Gamma((\alpha+1)i)} \right) A^{\ell} \ ; \ \forall p \geq \mu$$

$$= \sum_{\ell=0}^{p-1} \beta_{\alpha\ell}(t,p) A^{\ell} \ ; \ \forall p \geq \mu \qquad (3.4)$$

$\forall t \in \mathbf{R}_{0+}$, where, provided that $k-1 < \alpha(\in \mathbf{R}_+) \leq k(\in \mathbf{Z}_+)$, $\Gamma: \mathbf{R}_{0+} \to \mathbf{R}_+$ is the $\Gamma$-function of definition domain restricted to $\mathbf{R}_{0+}$ and the elements of the sets of functions $B_{\alpha j p} := \{\beta_{\alpha j\ell}: \mathbf{R}_{0+} \times \overline{p} \to \mathbf{R} : \ell \in \overline{p-1} \cup \{0\}\}$; $\forall j \in \overline{k-1} \cup \{0\}$ and $B_{\alpha p} := \{\beta_{\alpha\ell}: \mathbf{R}_{0+} \times \overline{p} \to \mathbf{R} : \ell \in \overline{p-1} \cup \{0\}\}$; $\forall p \geq \mu$ are defined for $t \in \mathbf{R}_{0+}$ as follows:

$$\beta_{\alpha j\ell}(t,p) = \left( \frac{t^{\alpha\ell}}{\Gamma(\alpha\ell+j+1)} + \sum_{i=p}^{\infty} \frac{\overline{a}_{\ell}(i) t^{\alpha i}}{\Gamma(\alpha i+j+1)} \right) ; \ \forall j \in \overline{k-1}, \ \forall \ell \in \overline{p-1} \cup \{0\}, \ \forall p \geq \mu \qquad (3.5)$$

$$\beta_{\alpha\ell}(t,p) = \left( \frac{t^{\alpha\ell}}{\Gamma((\alpha+1)\ell)} + \sum_{i=p}^{\infty} \frac{\overline{a}_{k}(\ell) t^{\alpha i}}{\Gamma((\alpha+1)i)} \right) ; \ \forall p \geq \mu \qquad (3.6)$$

**Remark 3.1**: Note that the homogeneous solution of the fractional differential system is given by:

$$x_{\alpha}(t) = \sum_{j=0}^{k-1} t^{j} \Phi_{\alpha j}(t) x_{0j} \qquad (3.7)$$

Then,

$$\left( {}^{C}D_{0+}^{\alpha} \left( \sum_{j=0}^{k-1} t^{j} \Phi_{\alpha j}(t) \right) \right)(t) = \sum_{j=0}^{k-1} \left( A t^{j} \Phi_{\alpha j}(t) \right) \qquad (3.8)$$

(see, for instance, [1] and [3]). On the other hand, if only one of the above point initial conditions $\varphi_{j}(0) = x_{j}(0) = x^{(j)}(0) = x_{0j}$ is nonzero for some arbitrary $j \in \overline{k-1} \cup \{0\}$ then, (3.8) is decomposable for each additive term resulting in:

$$x_{\alpha}(t) = t^{j} \Phi_{\alpha j}(t) x_{0j} \qquad (3.9)$$

Thus, $\Phi_{\alpha j}(t)$ is a fundamental matrix $\forall j \in \overline{k-1} \cup \{0\}$ which satisfies the fractional differential system of order $\alpha$:

$$\left( {}^C D_{0+}^\alpha t^j \Phi_{\alpha j} \right)(t) = A t^j \Phi_{\alpha j}(t)$$

$$= \left( {}^C D_{0+}^\alpha \left( \sum_{\ell=0}^{p-1} t^j \beta_{\alpha j \ell}(t,p) A^\ell \right) \right)(t) = \left( \sum_{\ell=0}^{p-1} t^j \beta_{\alpha j \ell}(t,p) A^{\ell+1} \right) \qquad (3.10)$$

; $\forall j \in \overline{k-1} \cup \{0\}$. Note that the fractional system becomes the standard one for $j = 0$ and $\alpha = 1$ so that $B_{10p} = B_p$ (defined in Lemma 2.4) and $\Phi(t) = e^{At}$. □

The following result of Section 2 extends the linear independence of the functions expanding $e^{At}$ to those appearing in the expansions of $\Phi_{\alpha j}(t)$; $\forall j \in \overline{k-1} \cup \{0\}$ and $\Phi_\alpha(t)$:

**Lemma 3.2**: All the sets of functions $B_{\alpha j p} := \{ \beta_{\alpha j \ell} : \mathbf{R}_{0+} \times \overline{p} \to \mathbf{R} : \ell \in \overline{p-1} \cup \{0\} \}$ and $B_{\alpha p} := \{ \beta_{\alpha \ell} : \mathbf{R}_{0+} \times \overline{p} \to \mathbf{R} : \ell \in \overline{p-1} \cup \{0\} \}$; $\forall p(\in \mathbf{Z}_+) \geq \mu$ are analytic and linearly independent on $\mathbf{R}_+$ for any given $\alpha \in \mathbf{R}_+$, $k \in \mathbf{Z}_+$ fulfilling $k-1 < \alpha \leq k$ and $\forall j \in \overline{k-1} \cup \{0\}$, $\forall p(\in \mathbf{Z}_+) \geq \mu$. Furthermore, the sets $B_{\alpha j p}$ satisfy a p-order linear time-invariant fractional differential system on $\mathbf{R}_+$, $\forall j \in \overline{k-1} \cup \{0\}$.

**Proof**: The sets of functions are analytic since from their defining formulas, it follows that they are infinitely differentiable on $\mathbf{R}_+$. Their linear independence is proven by contradiction. Since the rows of the fundamental matrices of solutions $\Phi_{\alpha j}(t)$; $\forall j \in \overline{k-1} \cup \{0\}$ of the Caputo fractional differential system of order $\alpha$ are linearly independent on $\mathbf{R}_{0+}$, it follows that $\Phi_{\alpha j}^T(t) \lambda = 0 \Leftrightarrow \lambda = 0$ for any $\lambda := (\lambda_1, \lambda_2, \cdots, \lambda_p)^T \in \mathbf{R}^p$; $\forall p \geq \mu$. But, if $B_{\alpha j p}$ is a linearly dependent set on $\mathbf{R}_{0+}$, since $\begin{bmatrix} I_n & A^T & \cdots & A^{p-1^T} \end{bmatrix}$ is a full row rank matrix, then it exists $\lambda \neq 0$ such that

$$t^{-j} \Phi_{\alpha j}^T(t) \lambda = \begin{bmatrix} I_n & A^T & \cdots & A^{p-1^T} \end{bmatrix} \begin{bmatrix} \beta_{\alpha j 0}(t,p) I_n \\ \beta_{\alpha j 1}(t,p) I_n \\ \vdots \\ \beta_{\alpha j, p-1}(t,p) I_n \end{bmatrix} \lambda = 0 \Leftrightarrow \begin{bmatrix} \beta_{\alpha j 0}(t,p) I_n \\ \beta_{\alpha j 1}(t,p) I_n \\ \vdots \\ \beta_{\alpha j, p-1}(t,p) I_n \end{bmatrix} \lambda = 0 \qquad (3.11)$$

; $\forall t \in \mathbf{R}_+$ which contradicts that $\lambda = 0$ in order that $\Phi_{\alpha j}(t)$ be a fundamental matrix, i.e. their rows are linearly independent functions. The functions in the set $B_{\alpha p} := \{ \beta_{\alpha \ell} : \mathbf{R}_{0+} \times \overline{p} \to \mathbf{R} : \ell \in \overline{p-1} \cup \{0\} \}$ are also linearly independent from the above considerations for the particular case $\alpha = k$, $j = k-1$ of $B_{\alpha j p}$. Note that the functions in the various sets are zero for t=0 so that linear independence is restricted to $\mathbf{R}_+$. On the other hand, the functions in the set $B_{\alpha j p}$ satisfy the differential system:

$$\left( {}^C D_{0+}^\alpha \left( \sum_{\ell=0}^{p-1} \beta_{\alpha j \ell}(t,p) A^\ell \right) \right)(t) = \left( \sum_{\ell=0}^{p-1} \beta_{\alpha j \ell}(t,p) A^{\ell+1} \right) \qquad (3.12)$$

; $\forall j \in \overline{k-1} \cup \{0\}$, subject to initial conditions $\beta_{\alpha j 0}(0,p) = 1/\Gamma(j+1) = 1/j!$, $\beta_{\alpha j \ell}(0,p) = 0$ ; $\forall j \in \overline{k-1}, \forall \ell \in \overline{p-1}$, $\forall p \geq \mu$. The functions in the sets $B_{\alpha p}$ satisfy initial conditions $\beta_{\alpha 0}(0,p) = 1/\Gamma(\alpha+1) = 1/(\alpha\Gamma(\alpha))$ $\forall j \in \overline{k-1}$, $\beta_{\alpha j \ell}(0,p) = \beta_{\alpha \ell}(0,p) = 0$; $\forall j \in \overline{k-1}, \forall \ell \in \overline{p-1}$, $\forall p \geq \mu$. □

## 4. Observability and controllability of linear fractional differential time-invariant systems

Consider the Caputo fractional differential system (3.1) with a measurable output defined by

$$y_\alpha(t) = C x_\alpha(t) \qquad (4.1)$$

for some $C \in \mathbf{R}^{z \times n}$ with $z(\in \mathbf{Z}_+) \leq n$. The following observability property is characterized.

**Definition 4.1**. The Caputo fractional differential system of order $\alpha$ is said to be observable in the observation time interval $[0,t]$ if $x_j(0) = x_{0j}$ ; $\forall j \in \overline{k-1} \cup \{0\}$ can be uniquely calculated from the measurable output $y_\alpha(t)$; $t \in [0,t]$ for some real interval $[0,t]$ of nonzero measure. □

**Theorem 4.2**. The Caputo fractional differential system (3.1), (4.1) of order $\alpha$ is observable in $[0,t]$ for any finite $t \in \mathbf{R}_+$ only if the standard dynamic system, i.e. that resulting when $\alpha = 1$, that is:

$$\dot{x}(t) = A x(t) + B u(t) \; ; \; y(t) = C x(t) \qquad (4.2)$$

is observable.

**Proof**: First, consider the homogeneous Caputo fractional differential system (3.1), i.e. $u \equiv 0$. One gets from (3.2) – (3.3):

$$y_\alpha(t) = \left( \sum_{j=0}^{k-1} t^j C \Phi_{\alpha j}(t) x_{0j} \right) = \left( \sum_{j=0}^{k-1} \sum_{\ell=0}^{p-1} t^j C \beta_{\alpha j \ell}(t,p) A^\ell x_{0j} \right)$$

$$= \overline{\beta}_\alpha^T(t,p) \, Block\,Diag\left[ \overbrace{\mathbf{Ob}(A,C,p) \vdots \cdots \vdots \mathbf{Ob}(A,C,p)}^{k} \right] \overline{x}_0 \qquad (4.3)$$

where $\overline{x}_0 := \left[ x_{00}^T, x_{01}^T, \ldots, x_{0,k-1}^T \right]^T$, and

$$\overline{\beta}_\alpha^T(t,p) := \left[ \beta_{\alpha 0 0}(t,p) I_s \; \cdots \; \beta_{\alpha 0, p-1}(t,p) I_s \; \cdots \; t^{k-1} \beta_{\alpha, k-1, 0}(t,p) I_s \cdots t^{k-1} \beta_{\alpha, k-1, p-1}(t,p) I_s \right]$$

(4.4)

$\forall t (\in \mathbf{R}_{0+})$ is a $s \times pks$ real matrix function ; $\forall p \geq \mu$ where $\mathbf{Ob}(A,C,\mu) := \begin{bmatrix} C \\ CA \\ CA^{\mu-1} \end{bmatrix}$. Since (4.2) is observable then

$$rank\ Block\ Diag\left[\overset{k}{\underbrace{Ob(A,C,p)\vdots\cdots\vdots Ob(A,C,p)}}\right] = kn \Leftrightarrow rank\ Ob(A,C,\mu) = rank\ Ob(A,C,p) = n$$

(4.5)

Furthermore, since the sets $B_{\alpha j p} := \{\beta_{\alpha j \ell} : R_{0+} \times \overline{p} \to R : \ell \in \overline{p-1} \cup \{0\}\}$ are linearly independent on $R_+$ from Lemma 3.2, the map from $R^{kn}$ to $R^s \times [0,t]$ ; $\forall t \in R_+$ from the initial conditions to the output-trajectory of the homogeneous system defined by (4.3) is injective. Thus, $\overline{x}_0$ can be uniquely calculated from $y: R^{kn} \to R^s \times [0,t]$. The sufficiency part has been proven. The necessity is obvious since if the rank in (4.4) is less than n, so that (4.2) is not observable, then the above mentioned map from $R^{kn}$ to $R^s \times [0,t]$ is not injective. The above considerations on (4.3)-(4.5) may be extended directly to the case that the admissible control u(t) is nonzero with its proof remaining valid, by replacing via (3.2) and (4.1):

$$y_\alpha(t) = Cx_\alpha(t) \to \overline{y}_\alpha(t) = y(t) - \int_0^t (t-\tau)^{\alpha-1} C\Phi_\alpha(t-\tau) Bu(\tau) d\tau \qquad \square$$

**Remark 4.3**. Note that Definition 4.1 and Theorem 4.2 characterize the observability property on any time interval of nonzero measure so that the property is independent of the time interval used for observation purposes. On the other hand, the property is independent of $\alpha$ because it implies the observability of the standard system (4.2) and it is implied by such an observability. $\square$

Remark 4.3 is formally enounced in the subsequent result:

**Corollary 4.4**. The observability property of the Caputo fractional differential system (3.1), (4.1) is independent of any real fractional order $\alpha$, with $k-1 < \alpha \leq k$ for any given $k \in Z_+$, and independent of the observation interval provided it is of nonzero measure. $\square$

The observability property can be tested as follows.

**Corollary 4.5**. The Caputo fractional differential system (3.1), (4.1) of any real order $\alpha$, with $k-1 < \alpha \leq k$ for any given $k \in Z_+$, is observable if and only if any of the three equivalent conditions holds:

1) $rank\left[C^T \vdots A^T C^T \vdots \cdots \vdots A^{T^{\mu-1}} C^T\right] = n$ (4.6)

2) $rank\left[\lambda I_n - A^T \vdots C^T\right] = n$; $\forall \lambda \in \sigma(A)$ (4.7)

where $\sigma(A)$ is the spectrum of A

3) $\int_0^t \Psi_\alpha^T(\tau) C^T C \Psi_\alpha(,\tau) d\tau \succ 0$, for any finite time interval $[0,t]$ of nonzero measure

$$\Psi_\alpha(t) := \left[\Phi_{\alpha 0}(t) \vdots t\Phi_{\alpha 1}(t) \vdots \ldots \vdots t^{k-1}\Phi_{\alpha,k-1}(t)\right]$$ (4.8)

**Proof**: The equivalence of Conditions 1-2 follows from the equivalence of the observability test for the standard system (i.e. $\alpha = k =1$) with the Popov- Belevitch – Hautus spectral observability test, [8-9], with the observability of any Caputo fractional system of order $\alpha$ with $k-1 < \alpha \leq k$ by using Corollary 4.4. The equivalence of Conditions 1-2 with Condition 3 follows from the fact that the observability of the fractional system of order $\alpha$ is uniform with respect to time (Theorem 4.2 or Corollary 4.4) and the fact that (3.7) into (4.1) with (4.8) yields for the unique (nontrivial) measurable output trajectory of the homogeneous fractional system of order $\alpha$ for each $\bar{x}_0 (\neq 0) \in \mathbf{R}^{nk}$:

$$0 < \int_0^t y_\alpha^T(\tau) y_\alpha(\tau) d\tau = \bar{x}_0^T \left( \int_0^t \Psi_\alpha^T(\tau) C^T C \Psi_\alpha(\tau) d\tau \right) \bar{x}_0 \tag{4.9}$$

for any time interval $[0,t]$ finite of nonzero measure. Thus, $\bar{x}_0(\neq 0) \in \mathbf{R}^{nk}$ is unique from (4.9), and then the system is observable, if and only if $\int_0^t \Psi_\alpha^T(\tau) C^T C \Psi_\alpha(\tau) d\tau \succ 0$. Thus, Condition 3 is equivalent to Conditions 1- 2. □

**Remark 4.6**. Note from Condition 3 of Corollary 4.5 that if the Caputo fractional system of order $\alpha$, with $k-1 < \alpha \leq k$, is observable then:

$$\int_0^t \tau^{2j} \Phi_{\alpha j}^T(\tau) C^T C \Phi_{\alpha j}(\tau) d\tau \succ 0 \Leftrightarrow \int_0^t \Phi_{\alpha j}^T(\tau) C^T C \Phi_{\alpha j}(\tau) d\tau \succ 0 ; \forall j \in \overline{k-1} \cup \{0\}$$

for any finite time interval $[0,t]$ of nonzero measure. This is easily seen by taking initial conditions such that $x_{0j} \neq 0$ for some $j \in \overline{k-1} \cup \{0\}$ while $x_{0i} = 0$, $\forall i(\neq j) \in \overline{k-1} \cup \{0\}$. □

**Remark 4.7**. Note that Theorem 4.2 is alternatively proven from the spectral observability test as follows. Since the rank in (4.7) can only be lost for $\lambda \in \sigma(A)$ then (4.7) is equivalent to:

$$rank \left[ \lambda I_n - A^T \vdots C^T \right] = n ; \forall \lambda \in \mathbf{C} \tag{4.10}$$

Using Agarwal's Laplace transform of (3.1), [2], the spectral observability of the Caputo fractional system of order $\alpha$ is lost if and only if $rank \left[ \lambda^\alpha I_n - A^T \vdots C^T \right] < n$ for some $\lambda \in \mathbf{C}$. Thus, note that:

$$rank \left[ \lambda I_n - A^T \vdots C^T \right] < n, \text{ some } \lambda \in \mathbf{C} \Rightarrow rank \left[ \lambda_1^\alpha I_n - A^T \vdots C^T \right] < n \text{ for } \lambda_1 = \lambda^{1/\alpha} \tag{4.11}$$

$$rank \left[ \lambda^\alpha I_n - A^T \vdots C^T \right] < n, \text{ some } \lambda \in \mathbf{C} \Rightarrow rank \left[ \lambda_1 I_n - A^T \vdots C^T \right] < n \text{ for } \lambda_1 = \lambda^\alpha \tag{4.12}$$

Then, the Caputo fractional system of order $\alpha$, with $k-1 < \alpha \leq k$, is (is not) spectrally observable (equivalently, observable) if and only if the standard system $\alpha = k = 1$ is (is not) spectrally observable. □

The controllability property is now discussed:

**Definition 4.8**. The Caputo fractional differential system of order $\alpha$ is said to be controllable in the time interval $[0,t]$ of nonzero finite measure if it exists an admissible control $u: [0,t] \to \mathbf{R}^m$ which steers the state –trajectory solution to any prescribed value $x^* = x_\alpha(t)$ for any given initial conditions. □

**Theorem 4.9**. The Caputo fractional differential system (3.1), (4.1) of any real order $\alpha$, with $k-1 < \alpha \leq k$ for any given $k \in \mathbf{Z}_+$, is controllable if and only if any of the three equivalent conditions holds:

1) $rank\left[ B \vdots AB \vdots \cdots \vdots A^{\mu-1} B \right] = n$  (4.13)

2) $rank\left[ \lambda I_n - A \vdots B \right] = n$; $\forall \lambda \in \sigma(A)$  (4.14)

where $\sigma(A)$ is the spectrum of A

3) $\int_0^t \Phi_\alpha(\tau) BB^T \Phi_\alpha(\tau) d\tau \succ 0$, for any finite time interval $[0,t]$ of nonzero measure

**Proof**: Assume the forced solution of (3.2) under zero initial conditions, namely, $x_j(0) = x_{j0} = 0$; $\forall j \in \overline{k-1} \cup \{0\}$ so that:

$$x_\alpha(t) = \int_0^t (t-\tau)^{\alpha-1} \Phi_\alpha(t-\tau) Bu(\tau) d\tau = \sum_{\ell=0}^{p-1} \left( \int_0^t (t-\tau)^{\alpha-1} \beta_{\alpha\ell}(t-\tau, p) A^\ell Bu(\tau) d\tau \right)$$

$$= \left( \int_0^t (t-\tau)^{\alpha-1} \left[ \beta_{\alpha 0}(t-\tau, p), \beta_{\alpha 1}(t-\tau, p), \cdots, \beta_{\alpha, p-1}(t-\tau, p) \right] u(\tau) d\tau \right) \left[ B^T \vdots BA^T \vdots \cdots \vdots B^T A^{p-1^T} \right]^T$$

$$= Co(A, B, p) \bar{\gamma}_{u\alpha}(t, p) \qquad (4.15)$$

; $\forall t \in \mathbf{R}_{0+}$, $\forall p \geq \mu$ after using (3.4) where $\gamma_{u\alpha i}(t,p) := \int_0^t (t-\tau)^{\alpha-1} \beta_{\alpha i}(t-\tau, p) u(\tau) d\tau$ are m-vector functions being dependent on the control; $\forall i \in \overline{p-1} \cup \{0\}$; $\forall t \in \mathbf{R}_{0+}$, and

$$Co(A, B, p) := \left[ B \vdots AB \vdots \cdots \vdots A^{p-1} B \right] \qquad (4.16a)$$

$$\bar{\gamma}_{u\alpha}(t,p) := \left[ \gamma_{u\alpha 0}^T(t,p), \gamma_{u\alpha 1}^T(t,p), \cdots, \gamma_{u\alpha, p-1}^T(t,p) \right]^T ; \forall t \in \mathbf{R}_{0+} \qquad (4.16b)$$

are, respectively, a $n \times pm$ real matrix and a real $pm$-vector function. Since the functions in $B_{\alpha p}$ are linearly independent and analytic on $\mathbf{R}_+$, $\forall p(\in \mathbf{Z}_+) \geq \mu$ for any given $\alpha \in \mathbf{R}_+$ (with $k-1 < \alpha \leq k$ - Lemma 3.2) then $Y_{u\alpha p} := \{ \gamma_{u\alpha i} : \mathbf{R}_{0+} \times \overline{p} \to \mathbf{R} : i \in \overline{p-1} \cup \{0\} \}$; $\forall p(\in \mathbf{Z}_+) \geq \mu$ is, by construction, a linearly independent set of analytic m-vector functions on a finite nonzero measure time interval $[0,t] \subset \mathbf{R}_+$ for any given $\alpha \in \mathbf{R}_+$, $k-1 < \alpha \leq k$ provided that the control is nonzero on some subset of nonzero measure of $[0,t]$. Thus, note from (4.15) that for arbitrary $\mathbf{R}^n \ni x^* = x_\alpha(t)$ for some linearly independent set $Y_{u\alpha p}$ on $[0,t]$ defined for some admissible nonzero control $u: [0,t] \to \mathbf{R}^m$ the Caputo fractional differential system of order $\alpha$ is controllable, independent of $\alpha$, if and only if

$$rank\left[Co(A,B,\mu) \vdots x^*\right] = rank\,Co(A,B,p) = rank\,Co(A,B,\mu) = n \qquad (4.17)$$

what follows from Rouché- Froebenius theorem from Linear Algebra. Then, the controllability is guaranteed by that of the standard system with $\alpha = k = 1$. The equivalence of (4.13) with the spectral controllability condition (4.14) follows from Theorem 4.8 and the well-known duality result and the state trajectory solution (3.2). The pair $(A,C)$ is observable in the sense that (4.6), equivalently (4.7), holds if its dual pair $(A^T, C^T)$ is controllable. The equivalence of Conditions 1-2 with Condition 3 is proven as follows from (4.14). Assume that the control is generated as $u(\tau) = B^T \Phi_\alpha(t-\tau)K$; $\tau \in [0,t]$ for some $K \in \mathbf{R}^m$. Then, the controllability constraint

$$x^* = x_\alpha(t) = \left(\int_0^t (t-\tau)^{\alpha-1} \Phi_\alpha(t-\tau)\,BB^T\Phi_\alpha(t-\tau)d\tau\right)K$$

is solvable for any prescribed $x^* \in \mathbf{R}^n$ if and only if

$$\int_0^t (t-\tau)^{\alpha-1}\Phi_\alpha(t-\tau)\,BB^T\Phi_\alpha(t-\tau)d\tau \succ 0 \Leftrightarrow \int_0^t \Phi_\alpha(t-\tau)\,BB^T\Phi_\alpha(t-\tau)d\tau \succ 0 \qquad (4.18)$$

with the control solution being

$$u(\tau) = B^T\Phi_\alpha(t-\tau)\left(\int_0^t (t-\tau)^{\alpha-1}\Phi_\alpha(t-\tau)\,BB^T\Phi_\alpha(t-\tau)d\tau\right)^{-1} x^*;\ \tau \in [0,t] \qquad (4.19)$$

Then, the Caputo fractional differential system of any order $\alpha$ is controllable independent of $\alpha$ on $[0,t]$ if $\int_0^t \Phi_\alpha(t-\tau)\,BB^T\Phi_\alpha(t-\tau)d\tau \succ 0$. Sufficiency has been proven. Necessity follows by contradiction. Assume that the system is controllable and $\int_0^t \Phi_\alpha(t-\tau)\,BB^T\Phi_\alpha(t-\tau)d\tau \succ 0$ fails for a given $[0,t]$ of nonzero finite measure. Then, the columns of the matrix function $\Phi_\alpha(t-\tau)\,Bu(\tau)$ are not linearly independent vector functions on $[0,t]$ for any admissible control $u:[0,t]\to \mathbf{R}^m$. Then, from (4.15)-(4.16) either $rank\left[B\vdots AB \vdots \cdots \vdots A^{\mu-1}B\right] < n$ which contradicts the controllability Condition 1 or the set $Y_{u\alpha p}$ is not linearly independent on $[0,t]$ for any admissible control what is impossible if $rank\left[B\vdots AB \vdots \cdots \vdots A^{\mu-1}B\right] = n$. Then, Condition 3 is equivalent to the equivalent Conditions 1-3. All the proof may be easily re-addressed for nonzero initial conditions by replacing :

$$x^* = x_\alpha(t) \to x^* = \left(x_\alpha(t) - \sum_{j=0}^{k-1} t^j \Phi_{\alpha j}(t)x_{0j}\right) \qquad \square$$

Note that, although the controllability property is independent of $\alpha$, the control (4.18) depends on $\alpha$. The matrices $Co(A,B,\mu) := \left[B\vdots AB \vdots \cdots \vdots A^{\mu-1}B\right]$ and $Ob(A,C,\mu) := \begin{bmatrix} C \\ CA \\ CA^{\mu-1} \end{bmatrix}$ are referred to as the controllability matrix of the pair (A, B) and the observability matrix of the pair (A,C), respectively, and have to be of rank n in order that the standard system by controllable (see Theorem 4.9, Eq.(4.12)), respectively, observable (see Corollary 4.5, Eq.4.6).

**Remark 4.10**. Note from Corollary 4.5 that the rank condition (4.6) requires the necessary condition $\mu \geq [n/s]$ which is then a necessary condition for observability. In the same way, the rank condition (4.13) in Theorem 4.9 requires $\mu \geq [n/m]$ which is then a necessary condition for controllability. □

**5. Observability and controllability under non-uniform sampling**

The results of the above section on controllability and observability are now extended for, in general, non- uniform sampling. It has proven in Section 4 that if the standard continuous– time system (i.e. $\alpha = k =1$) is controllable/observable then, any fractional system of real order $\alpha$ keeps the corresponding property from that of the standard system. It is now proven that the properties are still kept under non-uniform sampling for almost all choices of the sampling instants provided that their numbers are non-lees than the degree $\mu$ of the minimal polynomial of the matrix A. Such generic choices are possible from the following important property: the linearly independent matrix/vector functions which expand the evolution and control operators defining the state and output trajectory solutions of the Caputo fractional system of order $\alpha$ in polynomial functions of the matrix A in (3.1) are Chebyshev systems, [4-7].

**Theorem 5.1**. Assume that the standard system is observable with $\mu \geq [n/s]$ (controllable). Then, the sampled Caputo fractional system of order $\alpha$ is observable (controllable) for almost any choices of $p \geq \mu$ distinct sampling instants.

**Proof**: Note from (4.3) and (4.15) that the measurable output of the homogeneous system of the Caputo fractional differential system of order $\alpha$ and the transpose of its forced solution under zero initial conditions at any time instant $t_i \in \mathbf{R}_+$ are, respectively:

$$y_\alpha(t_i) = \bar{\beta}_\alpha^T(t_i, p) \text{ Block Diag}\left[\overset{k}{\overbrace{\mathbf{Ob}(A,C,p) \vdots \cdots \vdots \mathbf{Ob}(A,C,p)}}\right] \bar{x}_0 \qquad (5.1)$$

$$x_\alpha^T(t_i) = \bar{\gamma}_{u\alpha}^T(t_i, p) \mathbf{Co}(A,B,p)^T \qquad (5.2)$$

$\forall p(\in \mathbf{Z}_+) \geq \mu$, $\forall \alpha(\in \mathbf{R}_+) \in [k-1, k]$, with $k \in \mathbf{Z}_+$ where $\mathbf{R}^{s \times pks} \ni \bar{\beta}_\alpha^T(t_i, p) = \hat{\beta}_\alpha^T(t_i, p) I_s$, with $\hat{\beta}_\alpha^T(t_i, p) \in \mathbf{R}^{1 \times pk}$, and $\bar{\gamma}_\alpha^T(t_i, p) \in \mathbf{R}^{1 \times pm}$. Now, define

$$\bar{y}_{\alpha p}(t_1, t_2, \ldots, t_{pk}) := \left[y_\alpha^T(t_1, p), y_\alpha^T(t_2, p), \cdots, y_\alpha^T(t_{kp}, p)\right]^T \in \mathbf{R}^{pks} \qquad (5.3)$$

$$\bar{\beta}_{\alpha p}(t_1, t_2, \ldots, t_{pk}) := \left[\bar{\beta}_\alpha^T(t_1, p), \bar{\beta}_\alpha^T(t_2, p), \cdots, \bar{\beta}_\alpha^T(t_{kp}, p)\right] \in \mathbf{R}^{pks \times pks} \qquad (5.4)$$

for some set of sampling instants $\{t_1, t_2, \ldots, t_{pk}\}$; and

$$\bar{x}_{\alpha p}(t_1, t_2, \ldots, t_{pk}) := \left[x_\alpha(t_1, p), x_\alpha(t_2, p), \cdots, x_\alpha(t_{kp}, p)\right]^T \in \mathbf{R}^{pm \times n} \qquad (5.5)$$

$$\bar{\gamma}_{u\alpha p}(t_1, t_2, \ldots, t_{pm}) := \left[\gamma_{u\alpha}(t_1, p), \gamma_{u\alpha}(t_2, p), \cdots, \gamma_{u\alpha}(t_{pm}, p)\right]^T \in \mathbf{R}^{pm \times pm} \qquad (5.6)$$

for some set of sampling instants $\{t_1, t_2, ..., t_{pm}\}$. One gets from (5.3) via (5.1) and (5.4), and from (5.5) via (5.2) and (5.6) the following linear algebraic systems:

$$\bar{y}_{\alpha p}(t_1, t_2, ..., t_{pk}) = \bar{\beta}_{\alpha p}(t_1, t_2, ..., t_{pk}) \text{Block Diag}\left[\overset{k}{\overbrace{Ob(A,C,p) \vdots \cdots \vdots Ob(A,C,p)}}\right]\bar{x}_0 \tag{5.7}$$

$$\bar{x}_{\alpha p}^T(t_1, t_2, ..., t_{pk}) = \bar{y}_{u\alpha p}(t_1, t_2, ..., t_{pm}) Co(A,B,p)^T \tag{5.8}$$

Note that $\text{Block Diag}\left[\overset{k}{\overbrace{Ob(A,C,p) \vdots \cdots \vdots Ob(A,C,p)}}\right] \in R^{kps \times kn}$ is full rank equal to $kn$ if $p \geq \mu \geq [n/s]$ and $\text{rank } Ob(A,C,p) = \text{rank } Ob(A,C,\mu) = n$ (i.e. if the standard system with $\alpha = k = 1$ is observable). A unique solution exists to (5.7) if and only if the set of sampling instants $\{t_1, t_2, ..., t_{pk}\}$ is chosen so that the square real matrix $\bar{\beta}_{\alpha p}(t_1, t_2, ..., t_{pk})$ is non-singular provided that $\text{rank } Ob(A,C,p) = \text{rank } Ob(A,C,\mu) = n$ since then:

$$\text{rank}\left[\bar{y}_{\alpha p}(t_1, t_2, ..., t_{pk}) \vdots \bar{\beta}_{\alpha p}(t_1, t_2, ..., t_{pk}) \text{Block Diag}\left[\overset{k}{\overbrace{Ob(A,C,p) \vdots \cdots \vdots Ob(A,C,p)}}\right]\right]$$

$$\text{rank}\left[\bar{\beta}_{\alpha p}(t_1, t_2, ..., t_{pk}) \text{Block Diag}\left[\overset{k}{\overbrace{Ob(A,C,p) \vdots \cdots \vdots Ob(A,C,p)}}\right]\right] = kn \tag{5.9}$$

so that (5.7) is a compatible linear algebraic system (then, the map from $R^{kn}$ to $R^{pks}$ defined via (5.7) is injective) leading to a unique solution $\hat{\bar{x}}_0$ being a least-squares estimation of $\bar{x}_0$ given by:

$$\hat{\bar{x}}_0 = \left(\Omega^T(t_1, t_2, ..., t_{pk})\Omega(t_1, t_2, ..., t_{pk})\right)^{-1}\Omega^T(t_1, t_2, ..., t_{pk})\bar{y}_{\alpha p}(t_1, t_2, ..., t_{pk}) \tag{5.10}$$

where:

$$\Omega(t_1, t_2, ..., t_{pk}) := \bar{\beta}_{\alpha p}(t_1, t_2, ..., t_{pk}) \text{BlockDiag}\left[\overset{k}{\overbrace{Ob(A,C,p) \ldots \ldots Ob(A,C,p)}}\right] \tag{5.11}$$

from Rouché-Froebenius theorem for compatibility from Linear Algebra. Then, the observability property is preserved from the standard system for such a set of samples. Also, if $\text{rank } Co(A,C,p) = \text{rank } Co(A,C,\mu) = n$ if $p \geq \mu$, so that the standard system is controllable, then

$$\text{rank}\left[\bar{x}_{\alpha p}^T(t_1, t_2, ..., t_{pk}) \vdots \bar{y}_{u\alpha p}(t_1, t_2, ..., t_{pm}) Co(A,B,p)^T\right] = \text{rank}\left[\bar{y}_{u\alpha p}(t_1, t_2, ..., t_{pm}) Co(A,B,p)^T\right]$$

$$= \text{rank}\left[Co(A,B,p)\right] = \text{rank}\left[Co(A,B,\mu)\right] = n \tag{5.12}$$

in (5.8), if and only if the set of sampling instants $\{t_1, t_2, ..., t_{pm}\}$ is chosen so that the square real matrix $\bar{y}_{u\alpha p}(t_1, t_2, ..., t_{pm})$ is non-singular. Then, the controllability property is preserved from the standard system for such a set of samples. Since the linearly independent real matrix vector functions $\bar{\beta}_{\alpha i}(t,p)$ and linearly independent real vector functions $\gamma_{u\alpha i}^{T}(t,p)$; $\forall i \in \overline{p-1} \cup \{0\}$; $\forall p \geq \mu$ of domain $R_{0+}$ are also Chebyshev systems, [4-7] both $\bar{\beta}_{\alpha p}(t_1, t_2, ..., t_{pk})$ and $\bar{\gamma}_{u\alpha p}(t_1, t_2, ..., t_{pm})$ are non-singular matrices for almost all choices of the sampling instants. □

**Remark 5.2**. It suffices to choose the sampling instants in Theorem 5.1 as being mutually distinct and belonging to real intervals of the form $[\eta, \eta + \pi/\omega)$ where $\eta \in R_+$ is arbitrary and $\omega$ is an upper-bound of the maximum eigenfrequency of (3.1), i.e. the maximum absolute value of the imaginary parts of all the complex eigenvalues of A, if any. Otherwise, $\omega = +\infty$. This procedure guarantees that the corresponding coefficient matrices of Chebyshev systems are non-singular, [4-7]. Some choices are preferred if computation is required. For instance, it can be suitable to solve the linear algebraic system (5.10) so that the coefficient matrix (5.11) be as better conditioned as possible by the choice of the sampling instants. Then, the influence of numerical errors in the computation of the solution is minimized. It can be also suitable to take a number of sampling instants p large enough to make the estimated initial condition sufficiently close to the real one in the least-squares computation procedure. □

For algebraic solvability, it might be possible to reduce the number of sampling instants in (5.7)-(5.8) so that the coefficient matrix of the algebraic systems. The procedure is formalized as follows:

**Theorem 5.3**. The following properties hold:
**(i)** Assume that the standard dynamic system is observable. Then, the Caputo fractional linear system of order $\alpha$ is observable under non-uniform sampling as well from a set of existing sampling instants $ST := \{t_1, t_2, ..., t_{\hat{n}}\}$ with $\hat{n} := \max_{1 \leq i \leq s} n_i$ for a set of integers $n_i \in Z_{0+}$ fulfilling $\sum_{i=1}^{s} n_i = kn$ in such a way that each $i-th$ component of the measured output, i.e. $y_{\alpha i}(t_i)$, is observed at a set of sampling instants $ST_i := \{t_1, t_2, ..., t_{n_i}\} \subseteq ST$ (which is empty if $n_i = 0$); $\forall i \in \bar{s}$ which fulfils that the coefficient matrix $\hat{\beta}_{\alpha \hat{n}}(t_1, t_2, ..., t_{\hat{n}})$ of $\bar{\beta}$ - functions, defined below in the proof, of the associated linear algebraic system of grouped data $\bar{y}_{\alpha p}(t_1, t_2, ..., t_{\hat{n}})$ is non-singular.

**(ii)** Assume that the standard dynamic system is controllable. Then, the Caputo fractional linear system of order $\alpha$ is controllable under non-uniform sampling as well from a set of existing sampling instants $ST := \{t_1, t_2, ..., t_{\hat{n}}\}$ organized in such a way that for each control component a set of sampling instants $ST_i := \{t_1, t_2, ..., t_{n_i}\} \subseteq ST$ is used, which is empty if (which is empty if $n_i = 0$); $\forall i \in \bar{m}$ where $\hat{n} := \max_{1 \leq i \leq m} n_i$ for a set of integers $n_i \in Z_{0+}$ fulfilling $\sum_{i=1}^{m} n_i = n$ provided that the associated coefficient matrix is non-singular.

**Proof.** **(i)** For observability under non-uniform sampling with reduced number of samples, the square coefficient matrix of linearly independent matrix functions (5.4), i.e. $\bar{\beta}_{\alpha p}(t_1, t_2,...,t_{pk})$, of the algebraic linear system of equations (5.7) may be reduced to a minimum order being the less integer non less than $kn$ in order to estimate $\bar{x}_0 \in R^{kn}$ from a set of samples of $y_\alpha(t)$. The associated algebraic problem is formulated so that each component $y_{\alpha i}(t)$ is observed at a set of sampling instants $\{t_{i1}, t_{i2},..., t_{in_i}\}$ ($i \in \bar{s}$), which can be the empty set for some of the components $j \in \bar{s}$ (and then $n_j = 0$), with a choice of the s nonnegative integers $n_i$ ($i \in \bar{s}$) fulfilling the constraint $\sum_{i=1}^{s} n_i = kn$. Since the functions in the coefficient matrix are a Chebyshev system the resulting coefficient matrix may be generically constructed as being non-singular for a set of sampling instants $\{t_1, t_2,...,t_{\hat{n}}\}$ selected as $t_{ij}(j \in \bar{n}_i) = t_j$ if $j \in \bar{n}_i$ and $\bar{n}_i \leq j$; $\forall i \in \bar{s}$, $\forall j \in \hat{n}$ with $\hat{n} := \max_{1 \leq i \leq s} n_i$. Eq. (5.1) is replaced by the component-to-component corresponding set of equations:

$$y_{\alpha i}(t_i) = \bar{\beta}_{\alpha i}^T(t_j, j \in \bar{n}_i) G_i \bar{x}_0 ; \forall i \in \bar{s} \tag{5.13}$$

where

$$G_i := \text{Block Diag}\left[\mathbf{Ob}(A, C_i^T, \mu) \overset{k}{\vdots \cdots \vdots} \mathbf{Ob}(A, C_i^T, \mu)\right]; \forall i \in \bar{s} \tag{5.14}$$

so that the whole algebraic system (5.7) is replaced with

$$\bar{y}_{\alpha p}(t_1, t_2,...,t_{\hat{n}}) = \hat{\beta}_{\alpha\hat{n}}(t_1, t_2,...,t_{\hat{n}}) \text{Block Diag}\left[G_1 \overset{s}{\vdots \cdots \vdots} G_s\right] \bar{x}_0 \tag{5.15}$$

where

$$\hat{\beta}_{\alpha\hat{n}}(t_1, t_2,...,t_{\hat{n}}) := \text{Block Diag}\left[\bar{\beta}_{\alpha 1}^T(t_j, j \in \bar{n}_1) \vdots \cdots \vdots \bar{\beta}_{\alpha s}^T(t_j, j \in \bar{n}_s)\right] \tag{5.16}$$

Then, the mapping from $R^{kn} \to R^{kn}$ defined from (5.15) is injective if and only if $\text{rank } \mathbf{Ob}(A,C,\mu) = n$, requiring $s \geq [\mu/n]$, and, furthermore, the set of sampling instants $\{t_1, t_2,...,t_{\hat{n}}\}$ is chosen in such a way that $\text{Det } \hat{\beta}_{\alpha\hat{n}}(t_1, t_2,...,t_{\hat{n}}) \neq 0$. Infinitely many such sets of sampling instants always exist from the property of Chebyshev systems of the functions defining $\hat{\beta}_{\alpha\hat{n}}(t_1, t_2,...,t_{\hat{n}})$.

**(ii)** The proof is similar to that of property (i) and then its details are omitted. The guidelines for controllability under non-uniform sampling with a reduced number of sampling instants are as follows. One proceeds in a closed way the square coefficient matrix (5.6), i.e. $\bar{\gamma}_{u\alpha p}(t_1, t_2,...,t_{pm})$, of the algebraic linear system of equations (5.8) with a set of nonnegative integers $n_i$ ($i \in \bar{m}$), related to the input components, satisfying $\sum_{i=1}^{m} n_i = n$ and a set of sampling instants $\{t_1, t_2, ..., t_{\hat{n}}\}$ selected as $t_{ij}(j \in \bar{n}_i) = t_j$ if $j \in \bar{n}_i$ and $\bar{n}_i \leq j$; $\forall i \in \bar{m}$, $\forall j \in \hat{n}$ with $\hat{n} := \max_{1 \leq i \leq m} n_i$. □

## ACKNOWLEDGMENTS

The author is grateful to the Spanish Ministry of Education by its partial support of this work through Grant DPI2009-07197. He is also grateful to the Basque Government by its support through Grants GIC07143-IT-269-07and SAIOTEK S-PE08UN15.